\documentclass{article}
\usepackage{amssymb}
\usepackage{amsmath, amsthm}

\def\eq{\hspace*{-2.0mm}&=&\hspace*{-2.0mm}}

\newtheorem{corollary}{Corollary}[section]
\newtheorem{definition}{Definition}[section]

\newtheorem{lemma}{Lemma}[section]
\newtheorem{proposition}{Proposition}[section]
\newtheorem{remark}{Remark}

\title{The Bourguignon Laplacian and harmonic symmetric bilinear forms}

\author{
Vladimir Rovenski\footnote{Department of Mathematics, University of Haifa, Mount Carmel, Haifa, 31905, Israel,
E-mail address: vrovenski@univ.haifa.ac.il},
\ Sergey Stepanov
\footnote{Department of Mathematics, Russian Institute for Scientific
and Technical Information of the Russian Academy of Sciences,
20, Usievicha street, 125190 Moscow, Russia, E-mail address:  s.e.stepanov@mail.ru}
\ and \
Irina Tsyganok
\footnote{Department of Data Analysis and Financial Technologies,
Finance University, 49-55, Leningradsky Prospect, 125468 Moscow, Russia,
E-mail address: i.i.tsyganok@mail.ru}
}

\begin{document}

\date{}
\maketitle

\begin{abstract}
The theory of harmonic symmetric bilinear forms on a Riemannian manifold is an analogue of the theory of harmonic exterior differential forms on this manifold. To show this, we must consider every symmetric bilinear form on a Riemannian manifold as a one-form with values in the cotangent bundle of this manifold. In this case, there are the exterior differential and codifferential defined on the vector space of these differential one-forms. Then a symmetric bilinear form is said to be harmonic if it is closed and coclosed as a one-form with values in the cotangent bundle of a Riemannian manifold. In the present paper we prove that the kernel of the little known Bourguignon Laplacian is a finite-dimensional vector space of harmonic symmetric bilinear forms on a compact Riemannian manifold. We also prove that every harmonic symmetric bilinear form on a compact Riemannian manifold with non-negative sectional curvature is invariant under parallel translations. In addition, we investigate the spectral properties of the little studied Bourguignon Laplacian.
\end{abstract}



\section{\large Introduction}

First, we recall here some well-known facts of the theory of harmonic exterior differential forms on an $n$-dimensional Riemannian manifold $(M,g)$ (see, for example, \cite{11}). We write $d : C^{\infty}(\Lambda^{p} M)\to C^{\infty}(\Lambda^{p+1} M)$ for the familiar \textit{exterior differentiation operator} and the vector bundle $\Lambda^{p} M$ of exterior differentiation $p$-forms ($p=1,\, \ldots ,\, n-1$).
If $d\, \omega =0$, then the $p$-form $\omega \in C^{\infty}(\Lambda^{p} M)$ is said to be \textit{closed}.
The \textit{codifferentiation operator} $\delta : C^{\infty}(\Lambda^{p+1} M)\to C^{\infty}(\Lambda^{p} M)$ is defined as the formal adjoint of $d$. If $\delta \, \omega =0$, then the $\left(p+1\right)$-form
$\omega \in C^{\infty}(\Lambda^{p+1} M)$ is said to be \textit{coclosed}. Moreover, if $\omega \in {\rm Ker}\,d\bigcap {\rm Ker}\, \delta $, then the $p$-form $\omega $ is said to be \textit{harmonic}. Using the operators $d$ and $\delta $, one constructs the well-known Hodge-de Rham Laplacian $\Delta_{H} :=\delta \, d+d\, \delta $. Its kernel ${\rm Ker}\, \Delta_{H} $ is a finite-dimensional vector space over~the field of~real numbers of harmonic $p$-forms on a compact Riemannian manifold $(M,g)$. Moreover, every harmonic $p$-form on a compact Riemannian manifold $(M,g)$ with non-negative curvature operator $\bar{R} : \Lambda^{2} M\to {\Lambda }^{2} M$ is invariant under parallel translations. If the curvature operator $\bar{R}$ is non-negative everywhere and positive at some point of $(M,g)$, then every harmonic $p$-form is identically zero. In conclusion, we recall also that the spectral theory of the Hodge-de Rham Laplacian is well known (see, for example, \cite{SPR}).

Second, we will consider the theory of harmonic symmetric bilinear forms is an analogue of the theory of harmonic exterior differential forms
 (see, for example, \cite{9}). To show this, we must consider a symmetric bilinear form $\varphi \in C^{\infty}(S^{2} M)$ as a one-form with values in the cotangent bundle $T^{*} M$ on $M$. In particular, in accordance with the general theory, it is possible to determine an induced exterior differential
$d^{\nabla} : C^{\infty}(S^{2} M)\to C^{\infty}(\Lambda^{2} M\otimes T^{*} M)$ on the vector space of $T^{*} M$-valued differential one-forms. In particular, if $d^{\nabla } \, \varphi =0$ then the form
$\varphi \in C^{\infty}(S^{2} M)$ is said to be \textit{closed} \textit{bilinear form}. In this case, the $\varphi \in C^{\infty}(S^{2} M)$ is a \textit{Codazzi tensor}. We recall here that a symmetric bilinear form is called a Codazzi tensor (named after D.~Codazzi) if its covariant derivative is a symmetric tensor (see \cite{1}; \cite[p.~435]{2}). In addition, we will call the Codazzi tensor \textit{trivial} if it is a constant multiple of metric (see also \cite{1}). Next, let $\delta^{\nabla} : C^{\infty}(\Lambda^{2} M\otimes T^{*} M)\to C^{\infty} (S^{2} M)$ be the formal adjoin operator of the exterior differential $d^{\nabla } $
 (see \cite[p.~355]{2} and \cite{9}), then the form $\varphi \in C^{\infty}(S^{2} M)$ is said to be \textit{harmonic} if $\omega \in {\rm Ker}\, d^{\nabla} \bigcap {\rm Ker}\, \delta^{\nabla}$ (see \cite[p.~270]{9} and \cite[p.~350]{14}).

 Using the operators $d^{\nabla } $ and $\delta^{\nabla } $, J.-P. Bourguignon constructed the Laplacian $\Delta_{B} : =d^{\nabla} \delta^{\nabla } +\delta^{\nabla } \, d^{\nabla }$ (see \cite[p.~273]{9}).
 One can prove that its kernel ${\rm Ker}\, \Delta_{B} $ is a finite-dimensional vector space over~the field of~real numbers of harmonic symmetric bilinear forms on a compact Riemannian manifold $(M,g)$. In turn, we will prove that every harmonic symmetric bilinear form on a compact Riemannian manifold $(M,g)$ with non-negative sectional curvature is invariant under parallel translations. In addition, if the sectional curvature of $(M,g)$ is positive at some point of $(M,g)$, then every harmonic symmetric bilinear form is trivial.
 Moreover, in our paper we will investigate the spectral properties of the Bourguignon Laplacian $\Delta_{B}$.

\smallskip
 \textbf{Acknowledgments}. Our work was supported by the Foundation for Basic Research of the Russian Academy of Science, project 16-01-00756.

\section{The Bourguignon Laplacian and its spectral properties}

Let $(M,g)$ be a compact manifold (without boundary) then $L^{2}(M, g)$ denotes the usual Hilbert space of functions or tensors with the global product (resp. global norm)
\[
 \left\langle \, u,\, w\, \right\rangle =\int_{M}\, g\, \left(\, u,\, w\right) \, dv_{g}
 \quad({\rm resp}.\ \left\| \, u\, \right\|^{2} =\int_{M}\, g\, \left(\, u,\, u\right)\, dv_{g}  ),
\]
where the measure ${\rm dv}_{g}$ is the usual measure relative to $g$
(we will omit the term ${\rm dv}_{g} $). In this case, $H^{2}(M, g)$ denotes the usual Hilbert space of functions or tensors determined $(M, g)$ with two covariant derivatives in $L^{2}(M,g)$ and with the usual product and norm.

We will consider a symmetric bilinear form $\varphi \in C^{\infty}(S^{2} M)$ as a one-form with values in the cotangent bundle $T^{*} M$ on $M$. This bundle comes equipped with the Levi-Civita covariant derivative $\nabla $, thus there is an induced exterior differential $d^{\nabla} : C^{\infty}(S^{2} M)\to C^{\infty} (\Lambda^{2} M\otimes T^{*} M)$ on $T^{*} M$-valued differential one-forms such~as
\begin{equation} \label{GrindEQ__2_1_}
 d^{\nabla } \varphi \, \left(X,\, Y,\, Z\right) :=\left(\nabla_{X} \varphi \right)\, \left(Y,\, Z\right)-\left(\nabla_{Y} \varphi \right)\, \left(X,\, Z\right)
\end{equation}
for any tangent vector fields $X,\, Y,\, Z$ on $M$ and an arbitrary $\varphi \in C^{\infty}(S^{2} M)$.

\begin{remark}\rm
The theory on $T^{*} M$-valued differential one-forms can be found in papers and monographs from the following list \cite{FR}, \cite[p.~133--134; 355]{2}, \cite{9}, \cite[p.~338]{10}, \cite[p.~349--350]{14} and \cite{20}.
\end{remark}

J.-P. Bourguignon defined in \cite[p.~273]{9} the Laplacian $\Delta_{B} : C^{\infty}(S^{2} M)\to C^{\infty}(S^{2} M)$ by the formula $\Delta_{B} :=\delta^{\nabla } d^{\nabla } +d^{\nabla } \delta^{\nabla } $ where $\delta^{\nabla} : C^{\infty}(\Lambda^{2} M\otimes T^{*} M)\to C^{\infty}(S^{2} M)$ is the formal adjoin operator of the exterior differential $d^{\nabla } $. If $(M, g)$ is a compact Riemannian manifold then by direct computations yield we obtain the following integral formula:
\begin{equation} \label{GrindEQ__2_2_}
 \left\langle \Delta_{B} \, \varphi ,\, \varphi \right\rangle =\left\langle d^{\nabla } \varphi ,\, d^{\nabla } \varphi \right\rangle +\left\langle \, \delta^{\nabla } \, \varphi ,\, \, \delta^{\nabla } \, \varphi \right\rangle
\end{equation}
Based on this formula, we conclude that the \textit{Bourguignon Laplacian} $\Delta_{B} $ is a non-negative operator. On the other hand, by the general theorem on elliptic operators (see \cite[p.~464]{2}; \cite[p.~383]{10}) we have the orthogonal decomposition
\begin{equation} \label{GrindEQ__2_3_}
 C^{\infty}(S^{2} M) ={\rm Ker}\, \Delta_{B} \oplus {\rm Im}\, \Delta_{B}
\end{equation}
with respect to the global scalar product $\langle \, \cdot \; ,\; \cdot \, \rangle $. The first component of the right-hand side of \eqref{GrindEQ__2_3_} is the kernel ${\rm Ker}\, \Delta_{B} $ of the Bourguignon Laplacian $\Delta_{B} $. It is well known from \cite[p.~464]{2} that ${\rm Ker}\, \Delta_{B} $ is a finite-dimensional vector space over~the field of~real numbers. Next, an easy computation yields the \textit{Weitzenb\"{o}ck decomposition formula} (see also \cite{FR}; \cite[p.~355]{2}; \cite[p.~273]{9})
\begin{equation} \label{GrindEQ__2_4_}
 \Delta_{B} \, \varphi =\bar{\Delta }\, \varphi +B\,\varphi
\end{equation}
where $\bar{\Delta }=\nabla^{*} \nabla $ is the \textit{rough Laplacian} (see \cite[p.~52]{2}). The second component of the right-hand side of \eqref{GrindEQ__2_4_} is called the \textit{Weitzenb\"{o}ck curvature operator} for the Bourguignon Laplacian $\Delta_{B} $. Moreover, we known that it has the form $B\, \varphi : =\varphi \circ
{\rm Ric}-{\mathop{R}\limits^{\circ }} \, \varphi $ where $\circ $ is a composition of endomorphisms and ${\mathop{R}\limits^{\circ }} $ is the linear maps of $S^{2} M$ into itself such that (see \cite[p.~52]{2})
\begin{equation} \label{GrindEQ__2_5_}
 \big({\mathop{R}\limits^{\circ }} \, \varphi\big)(X,\, Y)
 =\sum\nolimits_{i=1,\ldots ,n}\; \varphi \, \left(R(X,\, E_{i} )\, Y,\, E_{i} \right)
\end{equation}
for the curvature tensor $R$ of $(M, g)$, for any $\varphi \in C^{\infty}(S^{2} M)$
and an arbitrary local orthonormal basis $E_1,\ldots, E_n$ of vector fields on $(M, g)$. 
In addition, $B\, \varphi$ direct verification yields that $B\, g=0$ and
\begin{equation}\label{GrindEQ__2_6_}
 {\rm trace}_{g} \, (B\, \varphi )=0.
\end{equation}
Then from \eqref{GrindEQ__2_4_} and \eqref{GrindEQ__2_6_} we obtain the identity
\begin{equation} \label{GrindEQ__2_7_}
 {\rm trace}_{g} \, \left(\, \Delta_{B} \, \varphi \right)
 =\bar{\Delta }\, \left({\rm trace}_{g} \varphi \right).
\end{equation}
Next, we will consider the spectral theory of the Bourguignon Laplacian $\Delta_{B} : C^{\infty}(S^{2} M)\to C^{\infty}(S^{2} M)$. Let $(M, g)$ be a compact Riemannian manifold and $\varphi$ be a non-zero eigentensor corresponding to the eigenvalue $\lambda $,
that is $\Delta_{B}\, \varphi =\lambda \, \varphi $ and $\lambda$ is a real nonnegative number. Then we can rewrite the formula $\Delta_{B} \, \varphi =\bar{\Delta}\,\varphi + B\,\varphi $ in the following form $\lambda \, \varphi =\bar{\Delta }\, \varphi +B\, \varphi $. In this case, from \eqref{GrindEQ__2_7_} we obtain
\begin{equation} \label{GrindEQ__2_8_}
 \bar{\Delta }\, \left(\, {\rm trace}_{g} \varphi \right)=\lambda \, \left(\, {\rm trace}_{g} \varphi \right),
\end{equation}
where $\bar{\Delta} : C^{\infty}(M)\to C^{\infty}(M)$ is the ordinary \textit{Laplacian} defined by the formula $\bar{\Delta }\, f=-\; {\rm div}\, \left(\, {\rm grad}\, f\right)$ for any $f\in C^{\infty}(M)$. In this case, the following equation holds:
\[
 \left\langle \, \bar{\Delta }\, \left(\, {\rm trace}_{g} \varphi \right), {\rm trace}_{g} \varphi \right\rangle =\left\langle \nabla \, {\rm trace}_{g} \varphi ,\; \nabla \, {\rm trace}_{g} \varphi \right\rangle .
\]
Therefore, $\bar{\Delta }\, \left(\, {\rm trace}_{g} \varphi \right)=0$ if and only if ${\rm trace}_{g} \varphi ={\rm const}$. In this case, if \eqref{GrindEQ__2_8_} holds for ${\rm trace}_{g} \varphi ={\rm const}$ and $\lambda \ne 0$, then ${\rm trace}_{g} \varphi $ must be zero. We proved the following lemma.

\begin{lemma}
Let $(M, g)$ be an $n$-dimensional $\left(\, n\ge 2\right)$ compact Riemannian manifold  and $\Delta_{B} \, \varphi =\lambda \, \varphi $ for the Bourguignon Laplacian $\Delta_{B} : C^{\infty}(S^{2} M)\to C^{\infty}(S^{2} M)$ and for a non-zero eigenvalue $\lambda $. If ${\rm trace}_{g} \varphi ={\rm const}$, then ${\rm trace}_{g} \varphi =0$. On the other hand, if ${\rm trace}_{g} \varphi $ is not constant, then ${\rm trace}_{g} \varphi $ is an eigenfunction of the Laplacian $\bar{\Delta} : C^{\infty}(M)\to C^{\infty}(M)$ such that $\bar{\Delta}\, \left(\, {\rm trace}_{g} \varphi \right)=\lambda \, \left(\, {\rm trace}_{g} \varphi \right)$.
\end{lemma}

Standard elliptic theory and the fact that the Laplacian $\bar{\Delta}: C^{\infty}(M)\to C^{\infty}(M)$ is a self-adjoint elliptic operator implies that the spectrum of $\bar{\Delta }$ consist of discrete eigenvalues $0=\bar{\lambda }_{0} <\bar{\lambda }_{1} <\bar{\lambda }_{2} <\ldots ,$ which satisfy the equation $\bar{\Delta }\, f_{i} =\bar{\lambda }_{i} \, f_{i} $ for $f_{i} \ne 0$ (see, for example, \cite{33}). Here we will focus on bounds on the first non-zero eigenvalue $\lambda_{1} $ imposed by the Riemannian geometry of $(M,g)$. The first lower bound for $\lambda_{1} $ was proved by Lichnerowicz \cite{L}. The \textit{Lichnerowicz theorem} is following: If $(M,g)$ is a compact Riemannian manifold of dimension $n\ge 2$, whose Ricci curvature satisfies the inequality ${\rm Ric}\, \ge (n-1)\, k>0$ for some constant $k>0$, then the first positive eigenvalue $\bar{\lambda}$ of the Laplacian $\bar{\Delta} : C^{\infty}(M)\to C^{\infty}(M)$ has the lower bound $\bar{\lambda }\ge n\, k$. Yang \cite{Y} generalized the previous result in the following form: Let $(M,g)$ be a compact Riemannian manifold of dimension $n\ge 2$ with ${\rm Ric}\ge (n-1)\, k\ge 0$ for some non-negative constant $k$ and diameter ${\rm D}(M)$, then the first positive eigenvalue $\bar{\lambda }$ of the Laplacian $\bar{\Delta} : C^{\infty}(M)\to  C^{\infty}(M)$ satisfies the lower bound $\bar{\lambda }\ge \frac{1}{4} \, (n-1)\, k +\pi^2/{D^{2} (M)} $.

On the other hand, by the spectral theory (see, for example, \cite{33}), the Bourguignon Laplacian
${\rm \Delta}_{B}$ has a discrete set of eigenvalues $\left\{\, \lambda_{a} \right\}$ forming a sequence $0=\lambda_{0} <\lambda_{1} <\lambda_{2} <\ldots$, and $\lambda_{a} \to +\infty $ as $a\to +\infty $.
Any eigenvalue of ${\rm \Delta }_{B}$ has finite multiplicity and an arbitrary $\lambda_{a}^{} $ for $a\ge 1$ is positive because $\Delta_{B} $ is a non-negative elliptic operator. Then as a corollary of the above Lichnerowicz and Yang theorems, we can formulate the following proposition.

\begin{proposition}
Let $(M, g)$ be a compact Riemannian manifold of dimension $n\ge 2$ and $\lambda $ be a positive eigenvalue of the Bourguignon Laplacian $\Delta_{B} : C^{\infty}(S^{2} M)\to C^{\infty}(S^{2} M)$ such that its corresponding eigentensor $\varphi \in C^{\infty}(S^{2} M)$
has a non-zero trace. If the Ricci curvature $(M, g)$ satisfies the inequality $Ric\, \ge (n-1)\, k>0$ for some positive constant $k$, then $\lambda $ has the lower bound $\lambda \ge n\, k$. On other hand, if $Ric\, \ge (n-1)\, k\ge 0$ for some non-negative constant $k$, then $\lambda $ satisfies the lower bound  $\lambda \ge \frac{1}{4} \, (n-1)\, k+\frac{\pi^{2} }{{\rm D}^{2}(M)}$, where ${\rm D}(M)$ is the diameter of $(M,g)$.
\end{proposition}

Next, we will consider the case of a positive eigenvalue $\lambda $ of the Bourguignon Laplacian $\Delta_{B} :  C^{\infty}(S^{2} M)\to C^{\infty}(S^{2} M)$ such that its eigentensor $\varphi $ is a traceless bilinear form.
In other words, $\varphi \in C^{\infty}(S_0^2 M)$ for the vector bundle of traceless symmetric bilinear forms $S_{0}^{2} M$.

Then, using \eqref{GrindEQ__2_4_}, we obtain the integral equality
\begin{equation}\label{GrindEQ__2_9_}
 \lambda \, \left\langle \varphi \, ,\, \; \varphi \right\rangle =\left\langle \, B\, \varphi \, ,\, \varphi \right\rangle +\, \left\langle \nabla \, \varphi \, ,\, \nabla \, \varphi \right\rangle .
\end{equation}
At the same time, by direct computations yield we obtain the following identity
\[
 g\, \left(\, B\, \varphi ,\; \varphi \right) = (1/2)\,g\left(K\, \varphi ,\; \varphi \right)
\]
where $K : =Ric\circ \varphi +\varphi \circ Ric-2{\mathop{\, R}\limits^{\circ }} \varphi $ is the \textit{Weitzenb\"{o}ck curvature operator} of the well known \textit{Lichnerowicz Laplacian} (see \cite[p.~54]{2}; \cite[p.~388]{10})
\begin{equation} \label{GrindEQ__2_10_}
 \Delta_{L} \, \varphi =\bar{\Delta }\, \varphi + K\, \varphi .
\end{equation}
In addition, direct verification yields that $K\, g=0$ and
\begin{equation} \label{GrindEQ__2_11_}
 {\rm trace}_{g} \, \left(\, K\, \varphi \, \right)=0.
\end{equation}
Let $\left\{\, e_{i} \right\}$ be an orthonormal basis of the tangent space $T_{x} M$ at an arbitrary point $x\in M$ such as $\varphi_{x} \left(\, e_{i} ,\, e_{j} \right)=\lambda_{i} \left(x\right)\, \delta_{ij} $ where $\delta_{ij}$ is the Kronecker symbol and ${\rm sec}\, \left(\, e_{i} \wedge \, e_{j} \right)$ be the sectional curvature in the direction of subspace $\pi \left(x\right)\subset T_{x} M$ for $\pi \left(x\right)={\rm span}\{\, e_{i} ,\, e_{j}\}$, then (see \cite[p.~388]{10})
\begin{equation} \label{GrindEQ__2_12_}
 g(K\, \varphi, \varphi)=\sum\nolimits_{\,i\ne j}\sec( e_{i} \wedge e_{j})(\varphi_{ii}-\varphi_{jj}).
\end{equation}
Let $S_{0}^{2} M$ be the vector bundle of traceless symmetric bilinear forms and $\Delta_{B} : C^{\infty}(S_{0}^{2} M)\to C^{\infty}(S_{0}^{2} M)$ be the Bourguignon Laplacian acting on the vector space of $C^{\infty}$-section of $S_{0}^{2} M$. If we denote by $K_{\rm min} $the minimum of the positive defined sectional curvature of $(M,g)$, i.e., ${\rm sec}\, \left(\sigma_{x} \right)\ge K_{\rm min} >0$ in all directions $\sigma_{\rm x} $ at each point $x\in M$, then from \eqref{GrindEQ__2_9_} we obtain the integral inequality
\begin{equation} \label{GrindEQ__2_13_}
 \lambda \, \left\langle \varphi ,\, \, \varphi \right\rangle  \ge \, \frac{1}{2} \, K_{{\it min}} \, \int_{M} \sum\nolimits_{i\; \ne j}\left(\, \varphi_{ii} -\varphi_{jj} \right)\,^{2}  \, dv_{g}  +\, \left\langle \nabla \, \varphi , \nabla \, \varphi \right\rangle \ge 0.
\end{equation}
for an arbitrary positive eigenvalue $\lambda $ corresponding to a non-zero eigentensor
$\varphi \in C^{\infty}(S_{0}^{2} )$ of $\Delta_{B}$. If the condition ${\rm trace}_{g} \, \varphi =\varphi_{11}^{} +\varphi_{22} +\ldots+\varphi_{n\, n}=0$ holds, then it is not difficult to prove the following equality
\[
 \left\| \, \varphi \, \right\|^{2} =\varphi_{11}^{2} +\varphi_{22}^{2} +\ldots +\varphi_{n\,n}^{2}
 =\frac{1}{n} \sum\nolimits_{i\,< j}\left(\varphi_{ii} -\varphi_{jj} \right)\,^{2}.
\]
because it equals to the following one:
\[
 \varphi_{11}^{2} +\varphi_{22}^{2} +\ldots +\varphi_{n\,n}^{2} = 
 -2 \sum\nolimits_{i\,< j} \varphi_{ii}\varphi_{jj},
\]
that is $(\varphi_{11}^{2} +\varphi_{22}^{2} +\ldots +\varphi_{n\,n}^{2})^2 = 0$.
In this case, from \eqref{GrindEQ__2_13_} one can obtain the integral inequality
\begin{equation} \label{GrindEQ__2_14_}
 (\lambda - n\, K_{\min })\int_{M}\| \varphi \|^{2} \, dv_{g}  \ge 0.
\end{equation}
Then from \eqref{GrindEQ__2_14_} we conclude that $\lambda \ge n\, K_{\min }$ for an arbitrary positive eigenvalue $\lambda $. In turn, if the first positive eigenvalue $\lambda =n\, K_{{\rm min}} $, then its corresponding traceless bilinear form $\varphi $ is invariant under parallel translation. In this case, if the holonomy of
$(M, g)$ is irreducible, then the tensor $\varphi $ must have the form $\varphi =\mu \cdot g$ for some constant $\mu $. But in our case, the identity holds ${\rm trace}_{g} \, \varphi =0$ and, consequently, we have $\mu =0$. Then the following statement holds.

\begin{proposition}
Let $(M,g)$ be an $n$-dimensional $(n\ge 2)$ compact Riemannian manifold and $\Delta_{B}: C^{\infty}(S_{0}^{2} M)\to C^{\infty}(S_{0}^{2} M)$ be the Bourguignon Laplacian acting on traceless symmetric bilinear forms. Then the first positive eigenvalue of $\Delta_{B} $ satisfies the lower bound $\lambda \ge n\, K_{\rm min} $ for the minimum $K_{\rm min} $ of the strictly positive sectional curvature of $(M, g)$. If the first positive eigenvalue $\lambda =n\, K_{\rm min}$, then the trace-free symmetric bilinear form $\varphi $ corresponding to $\lambda $ is invariant under parallel translation. In particular, if the holonomy of $(M,g)$ is irreducible, then this relation means that $\varphi \equiv 0$.
\end{proposition}

In particular, if $(M,g)$ is the standard sphere $\left(\, S^{n} ,\, g_{0} \right)$, then ${\rm sec}\,(X\wedge Y)=+1$ for orthonormal vector fields $X$and $Y$. In this case, the first positive eigenvalue $\lambda \ge n$. We can formulate the following corollary.

\begin{corollary}
Let $\left(\, S^{n} ,\, g_{0} \right)$ be the $n$-dimensional $(n\ge 2)$ standard sphere and $\Delta_{B} : C^{\infty}(S_{0}^{2} M)\to C^{\infty}(S_{0}^{2} M)$ be the Bourguignon Laplacian acting on traceless symmetric bilinear forms defined on $(S^{n}, g_{0})$.
Then the first positive eigenvalue of $\Delta_{B} $ satisfies the lower bound $\lambda \ge n$.
\end{corollary}

In the case of the standard sphere $\left(\, S^{n} ,\, g_{0} \right)$ we have $B\, \varphi: =\varphi \circ {\rm Ric}-{\mathop{R}\limits^{\circ }} \, \varphi = n\, \varphi $ and $K\varphi =2n\, \varphi $ for an arbitrary symmetric bilinear form $\varphi \in C^{\infty}(S_{0}^{2} M)$. Then we can write the equality $\Delta_{B} \, \varphi =\, \left(\, \mu - n\, \right)\, \varphi $ for an arbitrary positive eigenvalue $\mu $ of the Lichnerowicz Laplacian $\Delta_{L} $ and for some $\varphi \in C^{\infty}(S_{0}^{2} M)$ corresponding to $\mu $. It means that the eigenvalue $\lambda $ of $\Delta_{B} $, which corresponds to the same bilinear form $\varphi \in C^{\infty}(S_{0}^{2} M)$ is equal to $\lambda =\, \left(\, \mu - n\, \right)$. The converse is also true.

Consider the Lichnerowicz Laplacian $\Delta_{L} $ acting on traceless and divergence-free symmetric bilinear forms or, in other words, $TT$-\textit{tensors} defined on the standard sphere $\left(\, S^{n} ,\, g_{0} \right)$. In this case, we know from \cite{34} that the eigenvalues of $\Delta_{L} $ are given by the formula
$\mu_{a} =a(n-1+a)+2\,(n-1)$ for all $a\ge 2$, i.e.,
\[
 {\rm spec}\, \left(\, \Delta_{L} \left|{}_{TT} \right. \right)=\left\{\, a\, \left(\, n-1+a\, \right)+2\, \left(\, n-1\right)\, \, \left|\, a\ge 2\, \right. \right\}.
\]
Then we immediately obtain the spectrum of the $\Delta_{B} $ acting on the $TT$-tensors defined on the standard sphere $\left(\, S^{n} ,\, g_{0} \right)$:
\[{\rm spec}\, \left(\, \Delta_{B} \left|{}_{TT} \right. \right)=\left\{\, a\, \left(\, n-1+a\, \right)+\left(\, n-2\right)\, \, \left|\, a\ge 2\, \right. \right\}. \]
Based on this result, we can formulate the statement.

\begin{proposition}
The eigenvalues of the Bourguignon Laplacian $\Delta_{B} $ acting on the TT-tensors defined on the standard sphere $\left(\, S^{n} ,\, g_{0} \right)$ are given by the formula $\lambda_{a} =a\, \left(\, n-1+a\, \right)+\left(\, n-2\right)$ for $a\ge 2$.
\end{proposition}

\section{Harmonic symmetric bilinear forms and their vanishing theorems}

The formula \eqref{GrindEQ__2_1_} means that we take a symmetric bilinear form $\varphi \in C^{\infty}(S^{2} M)$ viewed as a one form with values in the tangent bundle. In this case, $\varphi\in C^{\infty}(S^{2} M)$ is a Codazzi tensor if and only if $d^{\nabla } \varphi =0$. Therefore, we can formulate the following obvious statement.

\begin{lemma}\label{L-2}
A symmetric bilinear form $\varphi \in C^{\infty}(S^{2} M)$ on a Riemannian manifold $(M, g)$ is a Codazzi tensor if and only if it is a closed one-form viewed as a one form with values in the tangent bundle $T^{*}M$ on $M$.
\end{lemma}

 J.-P. Bourguignon proved in \cite[p.~271]{9} that
\begin{equation}\label{GrindEQ__3_1_}
 \delta^{\nabla } \varphi =-\, d\, \left(\, {\rm trace}_{g} \varphi \, \right)
\end{equation}
for an arbitrary Codazzi tensor $\varphi \in C^{\infty}(S^{2} M)$. At the same time, he defined a \textit{harmonic symmetric bilinear form} in \cite[p.~270]{9}.

\begin{definition}\rm
A symmetric bilinear form $\varphi \in C^{\infty}(S^{2} M)$ on a Riemannian manifold $(M,g)$ is \textit{harmonic} if $\varphi \in {\rm Ker}\, d^{\nabla } \bigcap {\rm Ker}\, \delta^{\nabla } $.
\end{definition}

 Using Lemma~\ref{L-2} and equation \eqref{GrindEQ__3_1_}, this definition can be simplified slightly.

\begin{proposition}
A symmetric bilinear form $\varphi \in C^{\infty}(S^{2} M)$ on a Riemannian manifold $(M, g)$ is harmonic if and only if it is a Codazzi tensor with constant trace.
\end{proposition}

Based on the formula \eqref{GrindEQ__2_2_} and \eqref{GrindEQ__2_3_}, we conclude that the kernel of the \textit{Bourguignon Laplacian} $\Delta_{B} :=\delta^{\nabla } d^{\nabla } +d^{\nabla } \delta^{\nabla } $ has finite dimension and satisfies the condition ${\rm Ker}\, \Delta_{B} ={\rm Ker}\, d^\nabla  \bigcap {\rm Ker}\, \delta^{\nabla }$ on a compact Riemannian manifold $(M,g)$. Therefore, $\Delta_{B} $-harmonic bilinear forms are harmonic symmetric bilinear forms on a compact Riemannian manifold $(M,g)$. Therefore, we have the following.

\begin{proposition}
Let $(M,g)$ be an $n$-dimensional compact Riemannian manifold and
$\Delta_{B} : C^{\infty}(S^{2} M)\to C^{\infty}(S^{2} M)$ be the Bourguignon Laplacian.
Then the kernel of ${\rm \Delta}_{B} $ is the finite dimensional vector space of harmonic symmetric bilinear forms (or, in other words, Codazzi tensors with constant trace).
\end{proposition}

 J.-P. Bourguignon also proved in \cite[p.~281]{9} that a compact orientable Riemannian four-manifold admitting a non-trivial Codazzi tensor with constant trace must have \textit{signature} zero (see, for definition, \cite[p.~161]{2}). Then the follo\-wing corollary holds.

\begin{proposition}
Let $(M,g)$ be a compact orientable Riemannian four-dimensional mani\-fold. If the kernel of the Bourguignon Laplacian ${\rm \Delta}_{B} $ is non-trivial, then $(M,g)$ must have signature zero.
\end{proposition}

 Using the formula \eqref{GrindEQ__2_4_}, one can obtain the \textit{Bochner-Weitzenb\"{o}ck formula}
\begin{eqnarray} \label{GrindEQ__3_2_}
\nonumber
 (1/2)\Delta \| \varphi \|^{2} \eq -g(\bar{\Delta }\, \varphi, \varphi ) +\|\nabla \, \varphi\|^{2} \\
 \eq -g(\Delta_{B} \varphi,\, \varphi) +(1/2)\, g(K\, \varphi , \varphi) + \| \nabla \, \varphi\|^{2}
\end{eqnarray}
for an arbitrary $\varphi \in C^{\infty}(S^{2} M)$. Let $\varphi \in C^{\infty}(S^{2} M)$ be a harmonic form then \eqref{GrindEQ__3_2_} can be rewritten in the form (see also the formula \eqref{GrindEQ__2_12_})
\begin{equation}\label{GrindEQ__3_3_}
 \Delta\left\| \varphi \right\|^{2} =\sum\nolimits_{i\, \ne j}\sec\left(\, e_{i} \wedge \, e_{j} \right)\, \left(\, \varphi_{ii} -\varphi_{jj} \right) +2\left\| \nabla \, \varphi \right\|^{2} .
\end{equation}
We remind here that an arbitrary Codazzi tensor $\varphi $ commutes on $(M,g)$ with the Ricci tensor $Ric$ of $(M,g)$ at each point $x\in M$ (see \cite[p.~439]{2}). Therefore, the eigenvectors of an arbitrary Codazzi tensor $\varphi $ determine the principal directions of the Ricci tensor at each point $x\in M$ (see \cite[pp.~113--114]{21}). The converse is also true. Then taking into account of \eqref{GrindEQ__3_3_} and using the ``Hopf maximum principle", we will prove in the next paragraph that the following lemma holds.

\begin{lemma}\label{L-01}
Let $U$ be a connected open domain $U$ of a Riemannian manifold $(M,g)$ and $\varphi $ be a harmonic symmetric bilinear form defined at any point of $U$. If the sectional curvature ${\rm sec}\, \left(\, e_{i} \wedge \, e_{j} \right)\ge 0$ for all vectors of the orthonormal basis $\left\{\, e_{i} \right\}$ of $T_{x} M$ which is determined by the principal directions of the Ricci tensor ${\rm Ric}$ at an arbitrary point $x\in U$ and $\left\| \, \varphi \, \right\|^{2} $has a local maximum in the domain U, then $\left\| \, \varphi \, \right\|^{2} $ is a constant function and $\varphi $ is invariant under parallel translations in $U$. If, moreover, ${\rm sec}\, \left(\, e_{i} \wedge \, e_{j} \right)>0$ at some point $x\in U$, then $\varphi $ is trivial.
\end{lemma}

\proof Let us suppose that ${\rm sec}\, \left(\, e_{i} \wedge \, e_{j} \right)\ge 0$ in some connected open domain $U\subset M$ then $g\, \left(K\, \varphi ,\, \, \varphi \, \right)\ge 0$. If, moreover, there is a non-zero Codazzi tensor $\varphi$ given in $U\subset M$ then from \eqref{GrindEQ__3_3_} we conclude that $\Delta \, \left\| \, \varphi \, \right\|^{2} \ge 0$, i.e. $\left\| \, \varphi \, \right\|^{2} $ is a nonnegative subharmonic function in $U$. Let us suppose $\left\| \, \varphi \, \right\|^{2} $ has a local maximum at some point $x\in U$ then $\left\| \, \varphi \, \right\|^{2} $ is a constant function in $U\subset M$ according to the ``Hopf's maximum principle" (see \cite[p.~47]{22}). In this case, $\Delta \, \left\| \, \varphi \, \right\|^{2} =0$ and $\left\| \, \nabla \, \varphi \, \right\|^{2} =0$. In particular, the latter equation means that the form $\varphi $ is parallel.

 Let $\left\| \, \varphi \, \right\|^{2} =C$ for some constant $C$,  then from the equation \eqref{GrindEQ__3_3_} we obtain that $g\, \left(K\varphi ,\, \varphi \right)+2\, \left\| \, \nabla \, \varphi \, \right\|^{2} =0$. Since ${\rm sec}\, \left(\, e_{i} \wedge \, e_{j} \right)\ge 0$ it means that $g\, \left(K\, \varphi ,\, \varphi \right)=0$ and $\nabla \, \varphi =0$. If there is a point $x\in U$ such that ${\rm sec}\, \left(\, e_{i} \wedge \, e_{j} \right)>0$ then from \eqref{GrindEQ__3_3_} we come to the conclusion that $\lambda_{1}
 =\ldots =\lambda_{n} =\lambda $ which is equivalent to $\varphi = (1/n)\,g$, see \cite[p.~436]{2}.
 \hfill$\square$

\smallskip

If $(M,g)$ is a compact manifold and a harmonic symmetric bilinear form $\varphi $ is given in a global way on $(M,g)$ then due to the ``Bochner maximum principle" for compact manifold it follows the classical Berger-Ebin theorem (see \cite[p.~436]{2} and \cite[p.~388]{10}) which is a corollary of our Lemma~\ref{L-01}.

\begin{corollary}
Every harmonic symmetric bilinear form $\varphi \in C^{\infty}(S^{2} M)$ on a compact Riemannian manifold $(M,g)$ with nonnegative sectional curvature is invariant under parallel translations. Moreover, if ${\rm sec}>0$ at some point, then $\varphi \in C^{\infty}(S^{2} M)$ is trivial.
\end{corollary}

\begin{remark}\rm
It is well known that every parallel symmetric tensor field $\varphi \in C^{\infty}(S^{2} M)$ on a connected locally irreducible Riemannian manifold $(M,g)$ is proportional to $g$, i.e., $\varphi =\lambda\,g$ for some constant $\lambda $. Due to this the second parts of Corollary~3 can be reformulated in the following form: Moreover, if $(M,g)$ a connected locally irreducible Riemannian then an arbitrary harmonic symmetric bilinear form
$\varphi \in C^{\infty}(S^{2} M)$ is trivial.
\end{remark}

For example, let $(M,g)$ be a \textit{Riemannian symmetric space of compact type} that is a compact Riemannian manifold with non-negative sectional curvature and positive-definite Ricci tensor (see \cite[p.~256]{23}). Moreover, if a Riemannian symmetric space of compact type is a locally irreducibility Riemannian manifold $(M,g)$ then it is a compact Riemannian manifold with positive sectional curvature (see \cite{24}). Therefore, we can formulate the following corollary.

\begin{corollary}
Every harmonic symmetric bilinear form on a Riemannian symmetric manifold of compact type is invariant under parallel translations. If, in addition to the above mentioned the manifold is locally irreducible, then harmonic symmetric bilinear forms are trivial.
\end{corollary}

 The following theorem supplements the classical Berger-Ebin theorem (\cite[p.~436]{2} and \cite[p.~388]{10}) for the case of a complete noncompact Riemannian manifold.

\begin{proposition}
Let $(M,g)$ be a complete simply connected Riemannian manifold with nonnegative sectional curvature. Then there is no a non-zero harmonic symmetric bilinear form $\varphi \in C^{\infty}(S^{2} M)$ such that $\int_{M}\left\| \, \varphi \, \right\| \, d{\rm vol}_{g} <+\infty  $.
\end{proposition}

\proof
Let $(M,g)$ be a complete simply connected noncompact Riemannian manifold with nonnegative sectional curvature and $\varphi \in C^{\infty}(S^{2} M)$ be a globally defined non-zero harmonic symmetric bilinear form then $g(\, K\varphi ,\, \varphi)\ge 0$. Therefore, from \eqref{GrindEQ__3_3_} we obtain the inequality
\[
 \left\| \varphi \right\|\Delta \left\| \varphi \right\| =(1/2)\,g\left(K\varphi , \varphi \right) +
\left\| \nabla \, \varphi \right\|^{2} \ge 0.
\]
Then we conclude that $\left\| \varphi \right\|$ is a non-negative subharmonic function on a complete simply connected noncompact Riemannian manifold with nonnegative sectional curvature. In this case, if $\, \left\| \, \varphi \, \right\| $ is not identically zero, then it satisfies the condition $\int_{M}\left\| \varphi \right\| \, d{\rm vol}_{g} =\infty$ (see \cite{Wu}).
\hfill$\square$

\end{document}